\documentclass[a4paper,12pt
]{article}
\usepackage{amsmath}
\usepackage{amssymb}
\usepackage[dvipdfmx]{graphicx}
\usepackage{graphicx}
\usepackage{here}
\usepackage{mathpazo}

\newtheorem{thm}{Theorem}
\newtheorem{lem}[thm]{Lemma}

\newtheorem{ass}[thm]{Assertion}
\newtheorem{cor}[thm]{Corollary}

\newtheorem{prob}[thm]{Problem}

\newcommand{\R}{{\Bbb R}}
\newcommand{\N}{{\Bbb N}}

\newcommand{\IE}{{\it i.e.}, }
\newcommand{\EG}{{\it e.g.}, }
\newcommand{\RESP}[1]{[{\it resp.} {#1}]}

\newcommand{\ds}{\displaystyle}
\begin{document}
\begin{center}{
{\Large Derivatives of flat functions}
\vspace{10pt}
\\ 
Hiroki KODAMA, 
Kazuo MASUDA, 
and 
Yoshihiko MITSUMATSU
}
\end{center}
\begin{abstract}
We remark that 
there is no smooth function $f(x)$ 
on $[0, 1]$ which is flat at $0$ 
such that the derivative $f^{(n)}$ of any order $n\geq 0$ 
is positive on $(0,1]$. 
Moreover, the number of zeros of the $n$-th 
derivative $f^{(n)}$ grows to the infinity and the zeros accumulate 
to $0$ when $n \to \infty$.  
\end{abstract}

We consider smooth functions on the interval
$[0, 1]$ which are flat at the origin, 
namely of class $C^\infty$ and 
any derivative $f^{(n)}(x)$ converges to 0 when 
$x\to 0+0$.  Eventually it is equivalent to say that 
$f$ extends to the whole real line as a smooth function 
by defining $f(x)=0$ for $x<0$.  
In this short note we make a couple of remarks on the asymptotics 
of higher derivatives around the origin.   
The notion of higher derivatives was initiated in \cite{N}. 
For flat functions see, \EG \cite{GO}.

Among non-tirivial flat functions the most well-known 
might be the one which is defined as follows. 
$$
f(0)=0 \quad \mathrm{and}\quad   
f(x)=e^{-\frac{1}{x}}\,\,\, \mathrm{for}\,\,\, x>0
$$   
If we imagine its graph, 
of course it seems smooth enough,  
and it can be extended as constantly $0$ on 
$(-\infty, 0]$ as a smooth function 
on the real line $\R$. 
Its first derivative is positive on $(0, \infty)$, 
but the second derivative vanishes at 
$x=\frac{1}{2}=x_2$ 
and the third vanishes at $x_3=\frac{1-1/\sqrt 3}{2}<x_2$, 
and so on. 
That is, setting 
$x_n=\min\{x\,; \, f^{(n)}(x)=0,\, 
x>0\}$ 
for $n=2,\, 3,\, 4,\, \dots$, 
it is clear that $\{x_n\}_n$ is strictly decreasing, 
and in fact $\lim_{n\to\infty}x_n=0$. 
More over, if we fix any interval $[0, \alpha)$ ($\alpha>0$), 
$f^{(n)}(x)$ tends to behave more and more wildly  when $n\to\infty$ 
on the interval. 

Also, if we take $g_0(x)=f(x)(\sin(\frac{1}{x})+1)$ 
and 
$$g_n(x)=\int_0^x\int_0^{t_{n-1}}\cdots
\int_0^{t_1}
g_0(t_0)
dt_0\cdots dt_{n-2} dt_{n-1}\, ,
$$
then for $n=1,2,3,\cdots$, 
$g_n(x)$ is positive on $(0,\infty)$ and is flat 
at $x=0$, and apparently  $g_n^{(k)}(x)>0$ when $x>0$ 
for $0\leq k\leq n-1$ but there is no interval $(0, \alpha)$ 
on which $g_n^{(n)}(x)$ is positive. 


They seem to exhibit not particular for these examples 
but rather common or inevitable phenomena 
of higher derivatives of flat functions.

\begin{thm}{\rm \quad
There exists no smooth function $f(x)$ on $[0, 1]$ 
which is flat at $x=0$ and satisfies  
$f^{(n)}(x) >0$ on $(0, 1]$ for any $n\geq 0$. 
}
\end{thm}

This fact is refined as follows. 

\begin{thm}{\rm \quad
For a smooth function $f(x)$ on $[0, 1]$  
which is flat at $x=0$, put  
$Z(n)= \{x\in(0,1)\,\vert \, f^{(n)}(x)=0\}$ 
and $z(n)=\sharp Z(n)$ 
for $n\geq 0$. 
Then 
\vspace{-5pt}
$$ 
\lim_{n \to \infty}z(n)=\infty\, 
\vspace{-5pt}
$$
holds, where $\infty$ might be $\aleph_c$. 
}
\end{thm}

 

\begin{cor}{\rm \quad
1)\quad 
In general, 
$
\ds
\lim_{n\to\infty}
\inf Z(n) =0
$ . 
\\
2)\quad 
More strongly, for any $k>0$ there exist $N>0$ 
and 
$y^{(n)}(l)\in Z(n)$ for $n\geq N$ and $l= 1, \dots , k$ 
which are 
strictly increasing in $l$ and  
strictly decreasing in $n$, namely, 
satisfying 
\smallskip

for each fixed $n$,   
$\,\,y^{(n)}(l)<y^{(n)}(l+1)\,$ for $\,1\leq l \leq k -1$, 
\smallskip

for each fixed $l$, 
$\,\,y^{(n)}(l)>y^{(n+1)}(l)$.
\smallskip
\\
Moreover it satisfies for any $l$ 
$\ds \,\,\lim_{n\to\infty}y^{(n)}(l)=0$. 
}
\end{cor}
The accumulation of $Z(n)$ to $0$ ($n\to 0$) must be 
formulated in many more stronger statements. 
The above corollary is one of them.  
\medskip 
\\
These results imply that 
the set of germs of flat homeomorphisms 
of the half line $[0, 1)$ 
at $0$ 
is stable under the integration, 
while any element goes out of the set 
by a finite times of differentiation. 
\medskip
\\
\noindent
{\it Proof} of Corollary 3. \quad 
There is a zero of $f^{(n+1)}$ between 
two zeros of $f^{(n)}$. This simple argument,  
which will be used repeatedly, tells that 
once $Z(n)$ accumulates to $0$ for some $n$, 
so does $Z(k)$ for any $k\geq n$. 
Therefore in this case the proof is done.  
Otherwise, $0$ is always isolated from $Z(n)$ 
and then we can pick up the least element 
$y^{(n)}(1)\in Z(n)$.  
Now it is clear that $y^{(n+1)}(1)< y^{(n)}(1)$ for any 
$n$. 
Then, if $\ds \lim_{n\to\infty}y^{(n)}(1)=c>0$,  
$f\vert_{[0,c]}$ contradicts to Theorem 1. 
This proves 1).  

Now let us prove 2).  
Theorem 2 implies for any $k$ there is $N'$, 
$\sharp Z(n)\geq k$ for $n\geq N'$. Like in 1), 
once $0$ is accumulated by $Z(N)$, take any 
decreasing sequence $y^{(n)}(k)\in Z(n)$ for $n\geq N$, 
and then it is fairy easy to take $\{y^{(n)}(l)\}$ 
for other $l$'s 
so as to satisfy the conditions. 
Therefore we assume that 
$0$ is isolated from $Z(n)$ for  any $n\in \N$. 

Next, 
take $A(n)\subset Z(n)$ to be the set of points 
which is accumulated from above by points in $Z(n)$. 
Clearly this set has $\eta(n)=\min A(n)$ whenever 
$A(n)\ne\emptyset$. 
If $A(n)=\emptyset$, put $\eta(n)=1$.  
If $\eta(n)<1$, $f^{(n)}$ is flat at $\eta(n)$ and 
$\eta(n)\in A(n')$ for $n'\geq n$. 
Therefore the sequence $\{\eta(n)\}_n$ is weakly decreasing. 

In the case where $\ds c=\lim_{n\to\infty}\eta(n)>0$, 
applying Theorem 2 to $f\vert_{[0, c]}$, we can find $N$ 
such that $\sharp (Z(N)\cap(0,c)) \geq k$. 
Moreover, in this case, for any $n\geq N$ 
we can take the $k$ least zeros 
$0<y^{(n)}(1)<y^{(n)}(2)<\, \cdots \, < y^{(n)}(k)$ 
because there is no accumulation from above.  
Automatically $\{y^{(n)}(l)\}_n$ is strictly 
decreasing for each $l$.  
If $\ds \lim_{n\to\infty} y^{(n)}(k)=c'>0$, 
then again 
$f\vert_{[0,c']}$ contradicts to Theorem 2.  
Therefore this case is done. 

In the case where $\ds \lim_{n\to\infty}\eta(n)=0$, 
a similar argument in the case where $0$ is accumulated by 
some $Z(n)$ enable us to arrange $\{y^{(n)}(l)\}$ 
so as to satisfy the conditions. 
\hspace*{\fill}$\square$
\bigskip

\noindent
{\it Proof} of Theorem 1.  \quad 
The theorem is easily deduced from Lemma 4 by contradiction.   
Assume for some $\alpha\!>\!0$ that $f(x)$ is smooth on $[0, \alpha]$, is flat at $x=0$, 
and that its $n$-th derivative is positive on $(0, \alpha]$ for any 
$n\in \N$.  We adjust the function $f$ into 
$g(x)=f(\alpha)^{-1}f(\alpha x)$. 
Then $g(x)$ satisfies the condition of 
the lemma for any $n\in\N$. 
Therefore $g(x)\equiv 0$ on $[0, 1)$, 
and we obtain a contradiction. \hspace*{\fill} $\square$ 
\\
\begin{lem}
{\rm \quad Let $n$ be an integer and 
$g(x)$ be a function on $[0,1]$ 
of class $C^{n+1}$  
with the following properties. 
 \vspace{-5pt}
\begin{itemize}
\item[(1)]\quad $g^{(k)}(0)=0$ for $k=0,\, ...\, ,\, n$, and $g(1)=1$,  \vspace{-5pt}
\item[(2)]\quad $g^{(n+1)}(x)>0$ for $x>0$.  
 \vspace{-5pt}
\end{itemize}
Then \quad $g(x)<x^n$ \quad holds on $(0,1)$.  
}
\end{lem}

\noindent
{\it Proof} of Lemma 4. \quad 
It is enough to show that $g(x)/x^n$ is 
increasing on $[0, 1]$.  
As 
$\ds \frac{d}{dx}\left( \frac{g(x)}{x^n}\right)
=\frac{xg'(x) - ng(x)}{x^{n+1}}$, 
it is also sufficient to show that the numerator 
$xg'(x) - ng(x)$ is positive on $(0,1) $.  

Then because $(xg'(x) - ng(x))^{(n)}= xg^{(n+1)}(x)$  
 is positive on $(0,1]$ from our condition, 
we see successively that each $k$-th derivative 
 $(xg'(x) - ng(x))^{(k)}= xg^{(k+1)} - (n-k)g^{(k)}(x)$  
vanishes at $x=0$ and therefore is positive on $(0,1]$
for $k=n-1,\,n-2,\, \dots \, , 0$.   
This completes the proof. \hspace*{\fill} $\square$
\medskip
\\
A variant of this lemma is used to prove Theorem 2. 
\bigskip

\noindent
{\it Proof} of Theorem 2. \quad 
The key idea is not to look at $z(n)$ bt at the number 
$s(n)$ of the quasi-positive and quasi-negative intervals 
of $f^{(n)}$.   
 For a smooth (continuous) function $g$ on $[0,1]$ 
a connected component of the closure of 
$g^{-1}(0, \infty)$ \RESP{$g^{-1}(-\infty, 0)$} 
is called a {\it quasi-positive} \RESP{{\it quasi-negative}} 
interval. 
Such intervals are exactly maximal ones on which 
the primitive 
$\int g(x) dx$ of $g$  is strictly monotone. 
Let us define $s(g)\in \N\cup\{\infty\}$ 
to be the number of 
all the quasi-positive and quai-negative intervals of $g$. 
Then we put $s(n)=s(f^{(n)})$.  
\smallskip

If we have $s(n)=\infty$ for some $n$ ,  
$s(k)=\infty$ for $k\geq n$ as follows. 
If 
$Z(n)$ has interior points for some $n$, 
so does $Z(k)$ for $k\geq n$, 
If for some $n$ we have 
$z(n)=\infty$ and $\mathrm{int}\,Z(n)=\emptyset$, 
the complement $[0,1]\setminus Z(n)$ consists of 
infinitely many intervals and possibly of one half open interval.  
Each open interval contains an element in $Z(n+1)$.

Therefore, eliminating such cases, , 
we can assume $z(n)<\infty$ for any $n\in \N$. 
Consequently 
$s(n)<\infty$ ($\forall n\in \N$) holds as well.   
We want to prove 
$\ds \lim_{n\to\infty}s(n)=\infty$ 
under this assumption. 

Let $x^{(n)}(l)$ ($l\!\!=\!\!1,2,\dots,\,\, s(n)-1$, $n\!=\!\!0,1,2,\dots$) 
denotes the bigger end point 
of the $l$-th of quasi-positive/negative intervals 
for $f^{(n)}$, namely 
$[0,x^{(n)}(1)]$,  
$[x^{(n)}(1),x^{(n)}(2)]$,  
$\dots$\,, 
$[x^{(n)}(l-1),x^{(n)}(l)]$, 
$\dots$\,, 
$[x^{(n)}(s(n)-1),1]$ are the maximal intervals. 
Except for the final one, 
any quasi-positive 
\RESP{quasi-negative} interval contains 
a maximal \RESP{minimal} point in its interior. 
From this observation it is easy to see the 
following,   
among which 1) is a conclusion 
of Theorem 1, because it implies 
$z(n)\geq 1$ 
for some $n$ and then we have $s(n+1)\geq 2$.  
\begin{ass}{\rm
\quad 
1)\quad $s(n)\geq2$ for some $n$. 
\\
2)\quad $\{x^{(k)}(1)\}_k$ is strictly decreasing 
($k=n,\, n+1,\, n+2, \dots$)
for $n$ in 1).   
\\
3)\quad Also for any $m$ and 
$0<l<s(m)$, the sequence $\{x^{(k)}(l)\}_k$ 
is strictly 
\\
\qquad decreasing ($k=m,\, m+1,\, m+2, \dots$).  
\\
4)\quad $\{s(n)\}_n$ is weakly increasing, namely, 
$s(n)\leq s(n+1)$ for any $n$. 
}
\end{ass}
\smallskip

\noindent
Now let us procede by contradiction.   
We assume that $s(n)$ does not grow to $\infty$, \IE 
for some $N$,  $s(n)\equiv s(N)(=S)$ for any $n\geq N$. 
For fixed $l \in \{1, \cdots , S\}$, 
the quasi-positivity/negativity of the $l$-th 
interval is independent of $n\geq N$.  
Under the assumption we also see the following. 
\begin{ass}{\rm
\quad For $n\geq N$,  
\\
1)\quad
 $f^{(n)}$ is strictly monotone on the final interval 
$[x^{(n)}(S-1),1]$.  
\\
2)\quad
In particular $f^{(n)}(1)\ne 0$.  
More precisely, if the final interval is 
quasi-positive \RESP{quasi-negative} we have 
$f^{(n)}(1)> 0$ \RESP{$f^{(n)}(1)< 0$}. 
\\
3)\quad
The final intervals are increasing, namely, we have 
\\
\qquad
$x^{(N)}(S-1)>x^{(N+1)}(S-1)>\cdots
x^{(n)}(S-1)>x^{(n+1)}(S-1)>\cdots$. 
}
\end{ass}

By multiplying a non-zero constant to $f$, 
we assume that for any $n\geq N$ 
 $f^{(n)}$ is weakly 
increasing on the final interval and 
that $f^{(N)}(1)=1$.

\begin{lem}
{\rm \quad 
Under these assumptions, 
the following estimate holds.\vspace{-5pt}   
$$f^{(N)}(x)\leq x^p \quad \mathrm{on}\,\,\, 
[x^{(N)}(S-1),1]\, \quad 
\mathrm{for\,\,\, any}\,\,\,  p \in \N.
$$
}
\end{lem}

This lemma apparently implies 
$f^{(N)}\vert_{[x^{(N)}(S-1),1]} \equiv 0$ 
and contradicts to our assumption. 
This completes the proof of Theorem 2. 
\hspace*{\fill} $\square$
\bigskip

\noindent
{\it Proof} of Lemma 7. \quad 
We adjust the proof of Lemma 4 
in order to apply to $f^{(N)}$.  
Put $a_n=x^{(n)}(S-1)$ to simplify the notation. 

It is enough to show that 
for any $p\geq 0$, $f^{(N)}(x)\cdot x^{-p}$ 
is strictly increasing 
on $[a_N,1]$ because 
$f^{(N)}(x)\cdot x^{-p}\vert_{x=1}=1$.  
So it suffices to show   
$\ds \left(f^{(N)}(x)\cdot x^{-p}\right)'>0\,$,
namely, 
$\,xf^{(N+1)}(x)-pf^{(N)}(x)>0$ 
on $(a_N,1)$. 

For this purpose we prove inductively for 
$k=p,\, p-1,\, p-2,\, \cdots ,\, 0$  
\vspace{-8pt}
$$
\left(xf^{(N+1)}(x)-pf^{(N)}(x)\right)^{(k)}>0
\quad \mathrm{on} \quad (a_{N+k},1)\, .
\vspace{-8pt}
$$
For $k=p$, on $(a_{N+p+1},1)\,$ 
and in particular on $(a_{N+p},1)$, we have  
clearly \vspace{-8pt}
$$\left(xf^{(N+1)}(x)-pf^{(N)}(x)\right)^{(p)}=
xf^{(N+p+1)}(x)>0\, .$$ 
Then on each step,  
as 
$\,f^{(N+k+1)}(a_{N+k})>0\,$ and  
$\,f^{(N+k)}(a_{N+k})=0\,$, 
\vspace{-8pt}
$$\,
\left(xf^{(N+1)}(x)-pf^{(N)}(x)\right)^{(k)}
=
xf^{(N+1+k)}(x)-(p-k)f^{(N+k)}(x)
\vspace{-8pt}
$$ 
is positive at $\,x=a_{N+k}\,$.   
Because the inductive hypothesis implies 
its derivative is positive on 
$(a_{N+k},1)$ 
(and even on $(a_{N+k+1},1)$), 
the induction is completed. 
\hspace*{\fill} $\square$


\begin{prob}{\rm\quad
1)\quad For some smooth functions on $[0,1]$ which are flat 
at $0$, $\cup_{n=1}^\infty Z(n)$ seems to be dense in $[0,1]$. 
However, we do not see which kind of further properties 
as flat functions are essential for this phenomena, 
because it discusses points away from $0$. 
Verify this phenomena for certain $f$'s and 
explain the reason. 
\\
2)\quad 
Does there exist a smooth function on $[0,1]$ 
which is flat at $0$ 
such that $\lim_{n\to\infty}\max Z(n)=0$ ? 
Or 
how about flat at $0$ such that the derived set of 
$\cup_{n=1}^{\infty}Z(n)$ coincides with  $\{0\}$? 
It seems plausible that such functions do not exist, 
while we do not know how to prove it. 
}
\end{prob}

\mbox{}

\begin{flushright}
Hiroki KODAMA
\\
{\small \it
Graduate School of Mathematical Sciences, 
The University of Tokyo, 
\\
3-8-1 Komaba, Meguro-ku, 
Tokyo, 153-8914, Japan
\\
and
\\
Arithmer Inc.
\\
{\tt kodma@ms.u-tokyo.ac.jp}
\vspace{15pt}
}

Kazuo MASUDA
\\
{\small \it 3-40-15 Wakamiya, Nakano-ku, Tokyo 165-0033, Japan
\\
{\tt math21@maple.ocn.ne.jp}
\vspace{15pt}
}

Yoshihiko MITSUMATSU 
\\
{\small \it Department of Mathematics, 
Chuo University
\\
1-13-27 Kasuga Bunkyo-ku, 
Tokyo, 112-8551, 
Japan
\\
{\tt yoshi@math.chuo-u.ac.jp}
\\
}
\end{flushright}

\end{document}